\documentclass[11pt]{amsart}

\usepackage{a4wide,enumerate,color,graphicx}
\usepackage{amsmath,bm}
\usepackage{scrextend} 
\allowdisplaybreaks
\usepackage[normalem]{ulem}
\usepackage[dvipsnames]{xcolor}
\usepackage{setspace}
\usepackage{amsbsy}
\usepackage{hyperref}
\usepackage[overload]{empheq}

\newcommand{\R}{{\mathbb R}} 
\newcommand{\del}{\partial}
\newcommand{\diver}{\operatorname{div}} 

\newcommand{\dx}{\textnormal{d}x}
\newcommand{\dt}{\textnormal{d}t}

\newcommand{\dy}{\textnormal{d}y}

\newcommand{\T}{\mathbb{T}^3}
\newcommand{\D}{\mathbb{D}}

\numberwithin{equation}{section}
\hypersetup{
    colorlinks=true,
    linkcolor=blue,
    filecolor=blue,      
    urlcolor=blue,
}

\newtheorem{theorem}{Theorem}
\newtheorem{lemma}[theorem]{Lemma}

\newtheorem{definition}{Definition}

\begin{document}

	\title[Cahn--Hilliard/Navier--Stokes]{Energy identity for the incompressible Cahn--Hilliard/Navier--Stokes system \\
    with non--degenerate mobility}	

    \author{Stefanos Georgiadis}
    \address[Stefanos Georgiadis]{King Abdullah University of Science and Technology, CEMSE Division, Thuwal, Saudi Arabia, 23955-6900. \\
    and  Institute for Analysis and Scientific Computing, Vienna University of Technology, Wiedner Hauptstra\ss e 8--10, 1040 Wien, Austria}
    \email{stefanos.georgiadis@kaust.edu.sa} 

    \begin{abstract} 
        We consider the Cahn--Hilliard/Navier--Stokes system with non--degenerate mobility in the space--periodic case, describing the flow of two viscous immiscible and incompressible Newtonian fluids with matched densities. We identify sufficient conditions on the velocity field for weak solutions to satisfy an energy identity, improving previous results on the literature. 
    \end{abstract}

    \subjclass[2020]{35Q35}
    \keywords{Cahn--Hilliard/Navier--Stokes, weak solutions, energy identity, anomalous dissipation} 

	\maketitle

%%%%%%%%%%%%%%%%%%%%%%%%%%%%%%%%%%%%%%%%%%%%%%%%%%%%
%%%%%%%%%%%%%%%%%%%%%%%%%%%%%%%%%%%%%%%%%%%%%%%%%%%%

    \section{Introduction}

    Diffuse interface models have proven effective in explaining the evolution of two or more immiscible fluids. A key strength of these models is their capability to naturally depict topological changes, such as the merging or the separation of droplets. The foundational diffuse interface model in the case of two incompressible, viscous Newtonian fluids,  is referred to as Model H (cf. \cite{HoHa77}); it was derived also in \cite{GuPoVi96} in the framework of continuum mechanics and it was shown that it is consistent with the second law of thermodynamics. This model gives rise to the Navier--Stokes/Cahn--Hilliard system (see also \cite{AbGaGr12, LiNiSh20}):

    \noindent
    \begin{align}
        \label{NV} \del_tu+\diver(u\otimes u)-\diver(\nu(c)\D u)+\nabla p & = -\gamma\diver(\nabla c\otimes\nabla c), \\
        \label{incompr} \diver u & = 0, \\
        \label{CH} \del_tc+\diver(cu) & = \diver(m(c)\nabla\mu), \\
        \label{pot}\mu & = f(c) - \gamma\Delta c,
    \end{align} 

    \noindent
    where $u$ represents the volume--averaged velocity, $\D u = \frac{1}{2} (\nabla u + \nabla u^\top)$ stands for the viscous stress tensor, $p$ denotes the pressure, and $c \in [-1,1]$ is an order parameter linked to fluid concentration (e.g. concentration difference). Additionally, $\nu(c) > 0$ is the viscosity of the mixture, $\gamma > 0$ is a parameter related to the thickness of the interfacial region, $f=F'$, where $F$ represents a homogeneous free energy density, $\mu$ is the chemical potential, and $m = m(c) \geq 0$ the mobility coefficient.

    One of the fundamental assumptions of this model is that the densities of both components are the same, which restricts the applicability of the model to situations when density differences are negligible; for extensions of the above model in the case of different densities we refer to \cite{LoTr98, AbGaGr12, AbDeGa13}.

    To avoid complications with boundary terms, we restrict ourselves to the case of the 3--dimensional torus $\T := (\mathbb{R}/2\pi\mathbb{Z})^3$. Then, we consider equations \eqref{NV}--\eqref{pot} on $[0,T] \times \T$, where $T>0$, with the initial conditions

    \noindent
    \begin{equation}
        \label{IC} u(0,x)=u_0(x), \quad c(0,x)=c_0(x) \quad \textnormal{on } \T.
    \end{equation} 

    We place the following assumptions:

    \begin{itemize}
        \item[(A1)] $F \in C([-1,1]) \cap C^2((-1,1))$, with \[ \lim_{s\to-1}f(s) = -\infty, \quad \lim_{s\to1}f(s) = +\infty, \quad f'(s) \geq -\alpha, \textnormal{ for some $\alpha>0$}, \] 
        \item[(A2)] $\nu$ is a positive constant,
        \item[(A3)] $m$ is a positive constant.
    \end{itemize}

    \noindent
    We note that the most interesting cases in systems of Cahn--Hilliard type concern degenerate mobilities, something that our assumption (A3) excludes and we refer to \cite{AbDeGa13} for a treatment of the incompressible Cahn--Hilliard/Navier--Stokes system with degenerate mobility.  Nevertheless, the regularity provided by the existence theory is not enough and the fact that the mobility does not vanish is vital in our analysis.   

    Denoting by $L^2_\sigma$ the space of squared integrable solenoidal functions, we define a weak solution of \eqref{NV}--\eqref{IC} as follows:

    \begin{definition} \label{def}
        Let $u_0 \in L^2_\sigma(\T;\R^3)$ and $c_0 \in H^1(\T)$, with $|c_0| \leq 1$ almost everywhere on $\T$. Let assumptions (A1)--(A3) hold. A triplet $(u,c,\mu)$ with the regularity 
    
        \begin{align}
            & u \in C([0,T];L^2_\sigma(\T;\R^3)) \cap L^2(0,T;H^1(\T;\R^3)) \label{u} \\
            & c \in C([0,T];W^{1,q}(\T)), \textnormal{ for some $q>3$}, \quad \Delta c \in L^2(0,T;L^6(\T)) \label{c-der} \\
            & f(c) \in L^2(0,T;L^6(\T)), \quad |c| \leq 1 \textnormal{ a.e. on } \T \times (0,T) \label{c} \\
            & \nabla\mu \in L^2(0,T;L^2(\T;\R^3)) \label{mu}
        \end{align} 

        \noindent
        is called a weak solution of \eqref{NV}--\eqref{IC}, if it satisfies \eqref{NV}--\eqref{pot} in $\mathcal{D}^*((0,T)\times\T)$ and \eqref{IC} in $\mathcal{D}^*(\T)$ and the energy inequality 
        
        \begin{equation}
            \label{enin} \mathcal{E}(t) + \int_0^t\int_{\T} \left(\nu|\mathbb{D}u  |^2+m|\nabla\mu|^2\right) \dx\dt \leq \mathcal{E}(0)
        \end{equation} 

        \noindent
        holds almost everywhere in $(0,T)$, where the energy of the system is given by 
        
        \[ \mathcal{E}(t) = \int_{\T} \left(\frac{1}{2}|u|^2+\frac{\gamma}{2}|\nabla c|^2+F(c)\right) \dx. \]
    \end{definition}

    Weak solutions according to Definition \ref{def} assume a non--degenerate mobility and their existence has been established in \cite[Theorem 1]{Ab09}. In the case of a degenerate mobility coefficient, Definition \ref{def} no longer applies and we refer to \cite{AbDeGa13} for an appropriate definition.

    Our interest is in the energy study of system \eqref{NV}--\eqref{pot} under assumptions (A1)--(A3) and in particular, the goal is to identify conditions under which the energy inequality \eqref{enin} is satisfied as an identity. The 2--dimensional case has been treated in \cite{Ab09} and comes for free under no additional hypotheses. This is not the case for the 3--dimensional case. In \cite{LiNiSh20}, the authors show that in the case of a constant viscosity $\nu(c)=\nu$ and a constant mobility coefficient $m(c)=m$, weak solutions satisfy \eqref{enin} as an equality, under the additional assumptions \begin{equation}
    \label{LiNiSh1} f'\in C([-1,1])
    \end{equation} and \begin{equation}
    \label{LiNiSh2} \nabla u\in L^p(0,T;L^s(\T)),
    \end{equation} for specific exponents $p$ and $s$ (cf. \cite[Theorem 1.1]{LiNiSh20}).

    In the present work, we show that one can relax \eqref{LiNiSh2} and ask that the velocity field $u$ be only H\"older continuous in space. This idea originated from \cite{Be23} on the energy study for the Navier--Stokes equations. In particular, we prove the following theorem:
    
    \begin{theorem}\label{thm}
        Let $(u,c,\mu)$ be a weak solution according to Definition \ref{def}. Let assumptions (A1)--(A3) hold, along with \eqref{LiNiSh1}. Then, if \begin{equation}
            \label{vel-reg} u \in L^{\frac{2}{1+\alpha}+\delta}(0,T;C^\alpha(\T)), \textnormal{ for $\alpha > 1/3$ and $\delta>0$},
        \end{equation} the following energy equality holds true:
        \[ \mathcal{E}(t)+\int_0^T\int_{\T}\left(\nu|\mathbb{D}u|^2+m|\nabla\mu|^2\right)\dx\dt=\mathcal{E}(0), \quad \forall~t\geq0. \]
    \end{theorem}

\section{Proof of Theorem \ref{thm}}

Let $\rho \in C_0^\infty(\R^3)$ be a standard mollifier, i.e. $\rho\geq0$, with $\text{supp }\rho\subset B(0,1)\subset\R^3$, such that $\rho(x) = \rho(|x|)$ and $\int_{\R^{3}}\rho(x)\,\dx=1$. For $\epsilon\in(0,1]$, define
$\rho_{\epsilon}(x):=\epsilon^{-3}\rho(\epsilon^{-1}x)$. Then, for any function
$f\in L^{1}_{\textnormal{loc}}(\R^{3})$ we define by the convolution
\begin{equation*}
  f_{\epsilon}(x):=\int_{\R^{3}}\rho_{\epsilon}(x-y)f(y)\,\dy=\int_{\R^{3}}
  \rho_{\epsilon}(y)f(x-y)\,\dy. 
\end{equation*}

If $f\in L^{1}(\T)$, then $f\in L^{1}_{\textnormal{loc}}(\R^{3})$, and it turns out that $f_{\epsilon}\in C^{\infty}(\T)$ is $2\pi$--periodic along the $x_{j}$--direction, for
$j=1,2,3$. Moreover, if $f$ is a divergence--free vector field, then $f_{\epsilon}$ is a smooth divergence--free vector field.

In the following lemma we collect some properties that will be used throughout this work and for a proof we refer to \cite[Section 3]{Be23}.

\begin{lemma} \label{lemma}
If $u\in C^\alpha(\T)\cap L^1_{\textnormal{loc}}(\T)$, for some $\alpha\in(0,1)$, then
\begin{align}
    \label{eq:conv1}      
    & \max_{x\in\T}|u(x+y)-u(x)|\leq [u]_\alpha |y|^\alpha, \\
    \label{eq:conv2}  
    & \max_{x\in\T}|u-u_\epsilon|\leq [u]_\alpha \epsilon^\alpha,\\ 
    \label{eq:conv3}      
    & \max_{x\in\T}|\nabla u_{\epsilon}(x)|\leq C [u]_\alpha\epsilon^{\alpha-1}.
\end{align}
\end{lemma}

\noindent
In the sequel, we will also use the Constantin--E--Titi commutator from \cite{CoETi94}:
\begin{equation}
  \label{eq:commutator-Euler}
  (u\otimes u)_{\epsilon}=u_{\epsilon}\otimes u_{\epsilon}
  +r_{\epsilon}(u,u)-(u-u_{\epsilon})\otimes(u-u_{\epsilon}) ,  
\end{equation}
where
\begin{equation*}
  r_{\epsilon}(u,u):=\int \rho_{\epsilon}(y)[u(x-y)-u(x)]
  \otimes[u(x-y)-u(x)]\dy.
\end{equation*}

\noindent
Throughout this work, we set $\gamma = 1$ for convenience, since it does not affect our analysis.

We test \eqref{NV} against $\phi = (u_\epsilon)_\epsilon$. The argument would require a further smoothing with respect to time, which is standard (cf. \cite{BaTi18}) but not essential for the proof. Then, we get:

\begin{equation} \label{ns}
    \begin{split}
        \frac{1}{2}\|u_\epsilon(T)\|_{L^2(\T)}^2 & -\int_0^T\int_{\T}(u\otimes u)_\epsilon:\nabla u_\epsilon\dx\dt+\int_0^T\int_{\T}\nu|\nabla u_\epsilon|^2\dx\dt \\
        & = \frac{1}{2}\|u_\epsilon(0)\|_{L^2(\T)}^2+\int_0^T\int_{\T}(\nabla c\otimes\nabla c)_\epsilon:\nabla u_\epsilon\dx\dt
    \end{split}
\end{equation} and estimate the second term on the left--hand side and the second term on the right--hand side. For the former, we have: \[ \begin{split}
    I_1 & := \int_0^T\int_{\T}(u\otimes u)_\epsilon:\nabla u_\epsilon\dx\dt \\
    & = \int_0^T\int_{\T}[u_\epsilon\otimes u_\epsilon + r_\epsilon(u,u) - (u-u_\epsilon)\otimes(u-u_\epsilon)]:\nabla u_\epsilon\dx\dt \\
    & = \int_0^T\int_{\T}r_\epsilon(u,u):\nabla u_\epsilon\dx\dt - \int_0^T\int_{\T}(u-u_\epsilon)\otimes(u-u_\epsilon):\nabla u_\epsilon\dx\dt =: I_{11} + I_{12},
\end{split} \] since $u_\epsilon$ is divergence--free and thus \[ \int_0^T\int_{\T}u_\epsilon\otimes u_\epsilon:\nabla u_\epsilon\dx\dt = 0. \]

\noindent
Regarding the latter: \[ \begin{split}
    I_2 & := \int_0^T\int_{\T}(\nabla c\otimes\nabla c)_\epsilon:\nabla u_\epsilon\dx\dt \\
    & = - \int_0^T\int_{\T} (\nabla c_\epsilon\otimes\nabla c_\epsilon-(\nabla c\otimes\nabla c)_\epsilon):\nabla u_\epsilon\dx\dt-\int_0^T\int_{\T} u_\epsilon\cdot\diver(\nabla c_\epsilon\otimes\nabla c_\epsilon)\dx\dt \\
    & =: I_{21}+I_{22},
\end{split} \] where \[ I_{21} = -\int_0^T\int_{\T} (\Delta c_\epsilon\nabla c_\epsilon-(\Delta c\nabla c)_\epsilon+\nabla c_\epsilon\cdot\nabla(\nabla c)_\epsilon-\nabla c\cdot\nabla(\nabla c)_\epsilon)\cdot u_\epsilon\dx\dt, \] while due to \eqref{CH} \[ \begin{split}
    I_{22} & = - \int_0^T\int_{\T} (f(c)-\mu)_\epsilon u_\epsilon\cdot\nabla c_\epsilon\dx\dt - \int_0^T\int_{\T} u_\epsilon\cdot(\nabla c_\epsilon\cdot\nabla(\nabla c)_\epsilon)\dx\dt \\
    & = - \int_0^T\int_{\T} (f(c)-\mu)_\epsilon u_\epsilon\cdot\nabla c_\epsilon\dx\dt - \int_0^T\int_{\T} u_\epsilon\cdot\nabla\left(\frac{1}{2}|\nabla c_\epsilon|^2\right)\dx\dt \\
    & \overset{\eqref{incompr}}{=} - \int_0^T\int_{\T} (f(c)-\mu)_\epsilon u_\epsilon\cdot\nabla c_\epsilon\dx\dt.
    \end{split} \]

    Therefore, 
    \[ \begin{split}
    I_{22} & = - \int_0^T\int_{\T} f(c)_\epsilon u_\epsilon\cdot\nabla c_\epsilon\dx\dt + \int_0^T\int_{\T} \mu_\epsilon u_\epsilon\cdot\nabla c_\epsilon\dx\dt \\
    & = - \int_0^T\int_{\T} (f(c)_\epsilon-f(c_\epsilon)) u_\epsilon\cdot\nabla c_\epsilon\dx\dt - \int_0^T\int_{\T} f(c_\epsilon) u_\epsilon\cdot\nabla c_\epsilon\dx\dt \\
    & \phantom{xxx} + \int_0^T\int_{\T} \mu_\epsilon u_\epsilon\cdot\nabla c_\epsilon\dx\dt \\
    & = - \int_0^T\int_{\T} (f(c)_\epsilon-f(c_\epsilon)) u_\epsilon\cdot\nabla c_\epsilon\dx\dt - \int_0^T\int_{\T} u_\epsilon\cdot\nabla F(c_\epsilon)\dx\dt \\
    & \phantom{xxx} + \int_0^T\int_{\T} \mu_\epsilon u_\epsilon\cdot\nabla c_\epsilon\dx\dt \\
    & \overset{\eqref{incompr}}{=} - \int_0^T\int_{\T} (f(c)_\epsilon-f(c_\epsilon)) u_\epsilon\cdot\nabla c_\epsilon\dx\dt  + \int_0^T\int_{\T} \mu_\epsilon u_\epsilon\cdot\nabla c_\epsilon\dx\dt \\
    & := I_{221}+I_{222}.
\end{split} \] 

\noindent
Putting everything together, \eqref{ns} reads:

\begin{equation}\label{{eq1}}
    \begin{split}
        \frac{1}{2}\|u_\epsilon(T)\|_{L^2(\T)}^2 + I_{11} + I_{12} & + \int_0^T\int_{\T}\nu|\nabla u_\epsilon|^2\dx\dt \\
        & = \frac{1}{2}\|u_\epsilon(0)\|_{L^2(\T)}^2 + I_{21} + I_{221} + I_{222}.
    \end{split}
\end{equation}

Now, test \eqref{CH} against $(\mu_\epsilon)_\epsilon$, to obtain: \begin{equation} \label{ch}
\int_0^T\int_{\T}\mu_\epsilon\left(\del_tc_\epsilon+\diver(cu)_\epsilon-\diver(m\nabla\mu)_\epsilon\right)\dx\dt=0,
\end{equation} where \[ \begin{split}
    J_1 & := \int_0^T\int_{\T}\mu_\epsilon\del_tc_\epsilon\dx\dt \\
    & \overset{\eqref{pot}}{=} \int_0^T\int_{\T}(f(c)_\epsilon-\Delta c_\epsilon)\del_tc_\epsilon\dx\dt \\
    & = \int_0^T\int_{\T}(f(c)_\epsilon-f(c_\epsilon))\del_tc_\epsilon\dx\dt + \int_0^T\int_{\T}(f(c_\epsilon)-\Delta c_\epsilon)\del_tc_\epsilon\dx\dt \\
    & = \int_0^T\int_{\T}(f(c)_\epsilon-f(c_\epsilon))\del_tc_\epsilon\dx\dt + \int_0^T\int_{\T}\del_t\left(F(c_\epsilon)+\frac{1}{2}|\nabla c_\epsilon|^2\right)\dx\dt \\
    & =: J_{11} + J_{12}.
\end{split} \] 

\noindent
Moreover, 

\[ \begin{split}
    J_2 & := \int_0^T\int_{\T}\mu_\epsilon\diver(cu)_\epsilon\dx\dt \\
    & \overset{\eqref{incompr}}{=} \int_0^T\int_{\T}\nabla\mu_\epsilon\cdot(c_\epsilon u_\epsilon-(cu)_\epsilon)\dx\dt - \int_0^T\int_{\T}\mu_\epsilon u_\epsilon\cdot\nabla c_\epsilon\dx\dt \\
    & =: J_{21}+J_{22}
\end{split} \] and \[ J_3:= -\int_0^T\int_{\T}\mu_\epsilon\diver(m\nabla\mu)_\epsilon\dx\dt = \int_0^T\int_{\T} m|\nabla\mu_\epsilon|^2\dx\dt \] 

\noindent
and putting everything together, \eqref{ch} reads
\begin{equation} \label{{eq2}}
    J_{11} + J_{12} + J_{21} + J_{22} + J_3 = 0.
\end{equation}

\noindent
Then, we add up \eqref{{eq1}} and \eqref{{eq2}} and let $\epsilon \to 0$.

Notice that $I_{222}$ cancels out with $J_{22}$ and we start with $I_{12}$: \[ \begin{split}
    I_{12} & \leq \int_0^T\int_{\T} |u-u_\epsilon|^2|\nabla u_\epsilon| \dx\dt \\
    & = \int_0^T\int_{\T} |u-u_\epsilon|^{1+\alpha}|u-u_\epsilon|^{1-\alpha}|\nabla u_\epsilon|^\alpha |\nabla u_\epsilon|^{1-\alpha} \dx\dt, 
\end{split} \] and making use of \eqref{eq:conv1} and \eqref{eq:conv2} and after an application of H\"older's inequality

\[ \begin{split}
    I_{12} & \leq \int_0^T\int_{\T} |u-u_\epsilon|^{1+\alpha}\|u(t)\|_{C^\alpha(\T)}^{1-\alpha}\epsilon^{\alpha(1-\alpha)} C \|u\|_{C^\alpha(\T)}^\alpha \epsilon^{\alpha(\alpha-1)}|\nabla u_\epsilon|^{1-\alpha} \dx\dt \\
    & = C \int_0^T \|u\|_{C^\alpha(\T)} \left( \int_{\T} |u-u_\epsilon|^{1+\alpha}|\nabla u_\epsilon|^{1-\alpha} \dx \right) \dt \\
    & \leq C \int_0^T \|u\|_{C^\alpha(\T)} \|u-u_\epsilon\|_{L^2(\T)}^{1+\alpha}\|\nabla u_\epsilon\|_{L^2(\T)}^{1-\alpha} \dt.
\end{split} \] The fact that $u$ and thus $u_\epsilon$ are both in $L^\infty(0,T;L^2(\T))$, implies that $\|u-u_\epsilon\|_{L^2(\T)} \leq C$ for almost every $t\in(0,T)$, for some $C>0$ independent of $\epsilon$. Thus,

\[ \begin{split}
    I_{12} & \leq C \int_0^T \|u\|_{C^\alpha(\T)} \|\nabla u_\epsilon\|_{L^2(\T)}^{1-\alpha} \dt \\
    & \leq C \|u\|_{L^{\frac{2}{1+\alpha}}(0,T;C^\alpha(\T))}\|\nabla u\|_{L^2(0,T;L^2(\T))}^{1-\alpha},
\end{split} \] where in the second inequality we have used H\"older's inequality with exponents $\frac{1+\alpha}{2}$ and $\frac{1-\alpha}{2}$. Therefore, the integral is finite and since $(u-u_\epsilon)\otimes(u-u_\epsilon):\nabla u_\epsilon \to 0$ a.e. by properties of convolutions, we can use the dominated convergence theorem to show that $I_{12} \to 0$. The treatment of $I_{11}$ is identical, but instead of \eqref{eq:conv2}, we use \eqref{eq:conv1}.

Since $\nabla c\in L^2(0,T;W^{1,6}(\T))\cap L^\infty(0,T;L^3(\T))$, which implies $\Delta c\nabla c\in L^1(0,T;L^2(\T))$ and $u\in L^\infty(0,T;L^2(\T))$, we directly obtain that $I_{21}\to0$. 

Moreover, the terms $I_{221}, J_{11}$ and $J_{21}$ go to zero as $\epsilon\to0$. Indeed: \[ \begin{split}
    I_{221} & \leq \|\nabla c\|_{L^\infty(0,T;L^2(\T))}\|u\|_{L^2(0,T;L^3(\T))}\|f(c)_\epsilon-f(c_\epsilon)\|_{L^2(0,T;L^6(\T))}\to0.
\end{split} \] 

Regarding $J_{11}$, using equations \eqref{incompr} and \eqref{CH}, we get

\[ \begin{split}
    J_{11} & = - \int_0^T\int_{\T}(f(c)_\epsilon-f(c_\epsilon)) (u\cdot\nabla c)_\epsilon \dx\dt - m \int_0^T\int_{\T}\nabla(f(c)_\epsilon-f(c_\epsilon)) \cdot \nabla\mu_\epsilon \dx\dt \\
    & =: J_{111} + J_{112}.
\end{split} \] Clearly, $J_{111}\to0$ similar to $I_{221}$ and since $f'\in C([-1,1])$ and $\nabla\mu\in L^2(0,T;L^2(\T))$, $J_{112}\to0$. The fact that $\mu\in L^2(0,T;L^6(\T))$, $u\in L^\infty(0,T;L^2(\T))$ and $\nabla c\in L^\infty(0,T;L^3(\T))$ yield \[ \begin{split}
    J_{21} & \overset{\eqref{incompr}}{=} -\int_0^T\int_{\T}\mu_\epsilon(u_\epsilon\nabla c_\epsilon-(u\nabla c)_\epsilon)\dx\dt \\
    & \leq \|\mu\|_{L^2(0,T;L^6(\T))}\|u_\epsilon\nabla c_\epsilon-(u\nabla c)_\epsilon\|_{L^2(0,T;L^{6/5}(\T)}\to0.
\end{split} \] 

Finally, we show that \[ \begin{split} J_{12} & = \int_{\T}\left(F(c_\epsilon(T))+\frac{1}{2}|\nabla c_\epsilon(T)|^2\right)\dx -\int_{\T}\left(F(c_\epsilon^0)+\frac{1}{2}|\nabla c_\epsilon^0|^2\right)\dx \\
& \to \int_{\T}\left(F(c(T))+\frac{1}{2}|\nabla c(T)|^2\right)\dx -\int_{\T}\left(F(c^0)+\frac{1}{2}|\nabla c^0|^2\right)\dx.
\end{split} \]

\noindent
Indeed, for the terms involving $F$, we have that $c_\epsilon(x,T)\to c(x,T)$ almost everywhere on $\T$ and due to the continuity of $F$, $F(c_\epsilon(x,T)) \to F(c(x,T))$ almost everywhere on $\T$. Moreover, $|F(c_\epsilon(x,T))|\leq C$ for some positive constant $C$ independent of $\epsilon$, since $c_\epsilon$ is bounded and $F$ is continuous. The dominated convergence theorem can, then, be applied because constants are integrable on the torus. For the gradient terms, notice that due to \eqref{c-der}, $\nabla c$ is in $L^2(\T)$ uniformly in $t$ and thus $\nabla c_\epsilon(x,t) \to \nabla c(x,t)$ in $L^2(\T)$ for all $t\in[0,T]$ by convolution properties, in particular for $t=T$. Hence,
\[
    \begin{split}
        & \left| \int_{\T}|\nabla c_\epsilon(x,T)|^2\dx - \int_{\T}|\nabla c(x,T)|^2\dx \right| \\
        & \phantom{xxxxxxxxx} \leq \int_{\T}|\nabla c_\epsilon(x,T)-\nabla c(x,T)||\nabla c_\epsilon(x,T)+\nabla c(x,T)| \dx \\
        & \phantom{xxxxxxxxx} \leq \| \nabla c_\epsilon(x,T)-\nabla c(x,T) \|_{L^2(\T)} \| \nabla c_\epsilon(x,T)+\nabla c(x,T) \|_{L^2(\T)} \\
        & \phantom{xxxxxxxxx} \leq 2 \| \nabla c(x,T) \|_{L^2(\T)} \| \nabla c_\epsilon(x,T)-\nabla c(x,T) \|_{L^2(\T)} \to 0.
    \end{split}
\]
The kinetic energy terms and the terms at $t=0$ are treated in a similar way.

Putting everything together, we conclude the proof of the theorem.

%%%%%%%%%%%%%%%%%%%%%%%%%%%%%%%%%%%%%%%%%%%%%%%%%%%%
%%%%%%%%%%%%%%%%%%%%%%%%%%%%%%%%%%%%%%%%%%%%%%%%%%%%

\medskip

{\bf Acknowledgments.} The author would like to thank professor Luigi C. Berselli for pointing out the problem and for helpful conversations and the anonymous referee for his/her valuable comments.

\medskip

{\bf Funding.} Research partially supported by King Abdullah University of Science and Technology (KAUST) baseline funds. The author acknowledges partial support from the Austrian Science Fund (FWF), grants P33010 and F65. This work has received funding from the European Research Council (ERC) under the European Union's Horizon 2020 research and innovation programme, ERC Advanced Grant no. 101018153.

\medskip

{\bf Conflict of interest.} The author has no conflict of interest to report.

%%%%%%%%%%%%%%%%%%%%%%%%%%%%%%%%%%%%%%%%%%%%%%%%%%%%
%%%%%%%%%%%%%%%%%%%%%%%%%%%%%%%%%%%%%%%%%%%%%%%%%%%%

\end{document}